\documentclass[12pt,a4paper]{article}
\usepackage{float}
\usepackage{german}
\usepackage[german,english]{babel}
\usepackage[english]{babel}
\usepackage{epsfig}
\usepackage{latexsym}
\usepackage{epsf}
\usepackage {amssymb}
\usepackage {amsmath}
\usepackage{color}
\usepackage{setspace}
\newcommand{\leftsup}[2]{{\vphantom{#2}}_{#1}{#2}}
 
\usepackage{times}
\usepackage{bm}
\usepackage{natbib}

\newtheorem{lemma}{Lemma}
\newtheorem{theorem}{Theorem}
\newtheorem{corollary}{Corollary}

\bibpunct{(}{)}{;}{a}{,}{,}
\renewcommand{\cite}{\citep}
\newcommand{\biblist}{
\bibliographystyle{apalike}
\bibliography{EK}
}

\usepackage{epsfig}
\def\beqn{\begin{eqnarray*}}
\def\eeqn{\end{eqnarray*}}
\def\beq{\begin{eqnarray}}
\def\eeq{\end{eqnarray}}

\def\ss{\scriptsize}
\def\sett{\underline{\tau}_{\ss{K}}}
\def\sinc{\mbox{sinc}}
\def\i{{\tt{i}}}

\pdfoutput=1

\usepackage{rotating}

\begin{document}

\title{A unified framework for spline estimators}
\author{Katsiaryna Schwarz,\\Zurich Insurance Company Ltd, 
Mythenquai 2, \\CH-8022, Zurich, Switzerland,\\{Katsiaryna.schwarz@zurich.com}\\ \\
Tatyana Krivobokova\\
Institute for Mathematical Stochastics, \\Georg-August-Universit\"at G\"ottingen,
Goldschmidtstr. 7\\37077 G\"ottingen, Germany,\\ {tkrivob@gwdg.de}} 

\maketitle

\begin{abstract}
This article develops a unified framework
to study the asymptotic properties of all periodic spline-based
estimators, that is, of regression, penalized and smoothing splines. The explicit
form of the periodic Demmler--Reinsch basis in terms of
exponential splines allows the derivation of an expression for the
asymptotic equivalent
kernel on the real line for all spline estimators simultaneously. The corresponding
bandwidth, which drives the asymptotic behavior of spline estimators,
is shown to be a function of the number of knots and
the smoothing parameter. Strategies for the selection of the optimal bandwidth
 and other model parameters are discussed.
\end{abstract}


\section{Introduction}
\label{sec:intro}
Consider a nonparametric regression model for the data pairs
$(y_i,x_i)$
\beq
y_i=f(x_i)+\epsilon_i\quad (i=1,\ldots,N),
\label{eq:model}
\eeq
with the standard assumptions on the random errors, that is
$E(\epsilon_i)=0$ and $E(\epsilon_i\epsilon_j)=\sigma^2\delta_{ij}$
with $\sigma^2>0$ and $\delta_{ij}$ the Kronecker delta. Any
nonparametric linear
estimator of $f$ can be written as
$
\widehat{f}(x)={n}^{-1}\sum_{i=1}^nW(x,x_i)y_i.
$
For kernel estimators the weight function
$W(x,t)=h^{-1}K(x/h,t/h)$ for some bandwidth $h>0$ is
given explicitly \citep[][]{GasserMueller79} and is called the kernel function. This 
representation of $\widehat{f}$ allows a straightforward study of
its properties. 

Another class of nonparametric linear estimators is the spline-based
estimators: smoothing splines \citep[][]{Wahba75}, regression splines \citep[][]{AgarwalStudden80} and penalized splines \citep[][]{RuppertWandCarroll03}. Penalized 
splines combine projection onto a low-dimensional spline space with a
roughness penalty and circumvent certain practical
disadvantages of smoothing and regression splines.
For all spline-based estimators the exact
form of kernel $K(x,t)$ is unknown, but it can be sufficiently well
approximated for smoothing and regression splines. Such an
approximation is called the asymptotic equivalent kernel and is employed
to study the local asymptotic properties of spline estimators. \citet{CogburnDavis74}
obtained the asymptotic equivalent
kernel for smoothing splines on the real line using Fourier
techniques. \citet{MesserGoldstein93} and \citet{Agnan96}, see also \citet{Agnan04},
extended this kernel to the case of a bounded 
interval. \citet{EggermontLaRiccia06} refined these results. Other references on equivalent kernels for smoothing splines
are \citet{Silverman84}, \citet{Nychka95} and
\citet{AbramovichGrinshtein99}. Equivalent kernels for regression
splines have been derived
only on the real line in terms of B-splines of degree three in
\citet{HuangStudden93}.

The asymptotic properties of penalized spline estimators have received 
attention only recently. \citet{Claeskens-etal09} show that depending
on the number of knots taken, penalized splines have asymptotic
behaviour similar either to regression or to smoothing
splines.  \citet{Wangetal11} proved that for a certain number of knots the asymptotic
equivalent kernel of penalized spline estimators is asymptotically equivalent to that of the smoothing spline
estimators.  

In this article we aim to study all spline-based estimators in a unified framework. A new explicit expression for the Demmler--Reinsch basis for periodic
splines allows us to derive asymptotic equivalent kernels on the real line for all spline estimators that deliver  insights into the local
asymptotic behavior of spline estimators, depending on a certain
parameter that is a function of the number of knots and the smoothing
parameter.

\section{Model and equivalent kernels}
\label{sec:pensplines}
Consider the nonparametric regression model (\ref{eq:model}).
Let the data be equally spaced on the interval $[0,1]$,
i.e., $x_i=i/N$ ($i=1,\ldots,N$). The unknown regression function $f$ is assumed
to be sufficiently smooth. More
precisely, $f\in
{\mathcal{H}}_{p+1}[0,1]=\{f:\;f\in
{\mathcal{C}}^{p},\;\int_0^1 \{f^{(p+1)}(x)\}^2dx<\infty\}$. To
estimate $f\in{\mathcal{H}}_{p+1}[0,1]$, first define
a partition of $[0,1]$ into $K\leq N$ equidistant
intervals
$\underline{\tau}_{\ss{K}}=\{0=\tau_0<\tau_1<\cdots<\tau_{\ss{K-1}}<\tau_{\ss{K}}=1\}$
with $\tau_i=i/K$ ($i=0,\ldots,K$). 
The spline space ${\mathcal{S}}(p;\sett)$ of degree $p>0$
based on $\sett$ consists of functions $s\in
{\mathcal{C}}^{p-1}[0,1]$, such that $s$ is a degree $p$ polynomial on each
$[\tau_{i},\tau_{i+1})$ ($i=0,\ldots,K-1$). 
The spline estimator
$\widehat{f}$ is the solution to
\beq
\min_{s\in
  {\mathcal{S}}(p;\sett)}\left[\frac{1}{N}\sum_{i=1}^{N}\{Y_i-s(x_i)\}^2+{\lambda}\int_0^1\{s^{(q)}(x)\}^2dx\right],\;\;\lambda\geq
0,\;0<q\leq p.
\label{eq:min_problem}
\eeq
For
$K=N$ and $p=2q-1$ the solution to (\ref{eq:min_problem}) is the
smoothing spline estimator. If
$\lambda=0$ and $K\ll N$, (\ref{eq:min_problem}) yields the
regression spline estimator, and a
general estimator with $K<N$, $p+1 > q>0$ and $\lambda>0$ is the
penalized spline estimator.   
A solution to (\ref{eq:min_problem}) can be written as 
$
\widehat{f}(x)={N}^{-1}\sum_{i=1}^n{{W}}^{[0,1]}\left( x,x_i \right) Y_i,
$
where the bandwidth $h=h(n)>0$ decays to zero such that
$nh(n)\rightarrow\infty$. The weight function $W^{[0,1]}(x,t)$ depends on the
observations $x_i$ and on $N$ and therefore is called the effective kernel. Let ${\cal{W}}^{[0,1]}(x,t)$ denote
an approximation to $W^{[0,1]}(x,t)$, which is independent of $x_i$, $N$ and such that ${E}\left\{\sup_{x,q}\left|\widehat{f}(x)-N^{-1}\sum_{i=1}^N{\cal{W}}^{[0,1]}(x,x_i)Y_i\right|\right\}$
is negligible compared to the bias of $\widehat{f}$. Then,
${\cal{W}}^{[0,1]}(x,t)$ is called the asymptotic equivalent kernel for $\widehat{f}$.
Typically, ${\cal{W}}^{[0,1]}(x,t)$ is found
as a sum ${\cal{W}}^{[0,1]}(x,t)={\cal{W}}(x,t)+{\cal{W}}^b(x,t)$. Here,
${\cal{W}}(x,t)$ corresponds to the asymptotic equivalent kernel on
the real line and ${\cal{W}}^b(x,t)$ is the asymptotic boundary
kernel, which decays exponentially away from the boundaries. In
particular, $\widehat{f}(x)\approx
N^{-1}\sum_{i=1}^N{\cal{W}}(x,x_i)Y_i$ for $x$ away from the
boundaries.

The available results on equivalent kernels for spline estimators can
be summarised as follows. 
For smoothing splines ${\cal{W}}^{[0,1]}_{ss}(x,t)$ has been approximated
by the Green's function for the Euler equations
\begin{equation}
\begin{aligned}
f(x)+\lambda (-1)^q f^{(2q)}(x)&=g(x),\;\;x\in(0,1),\;\lambda>0,\\
f^{(j)}(0)=f^{(j)}(1)&=0\;\;(j=q,\ldots,2q-1).
\end{aligned}
\label{eq:Euler}
\end{equation}
The solution to (\ref{eq:Euler}) is found in
two steps. First, a fundamental solution is obtained as the
Green's function of Euler equation (\ref{eq:Euler}) with 
$\lim_{|x|\rightarrow\infty}f^{(j)}(x)=0$ ($j=q,\ldots,2q-1$).
This Green's function is the
asymptotic 
equivalent kernel ${\cal{W}}_{ss}(\cdot,t)$ on the real line  and is
known explicitly for each $q$. In particular, scaling with $h=\lambda^{1/(2q)}>0$ gives
\beq
h{\cal{W}}_{ss}\left(hx,ht\right)=\mathcal{K}_{{ss}}(x,t)=\sum_{j=0}^{q-1}\frac{\i\exp\left[\i\left|x-t\right|\exp\left\{ \pi\i\left(2j+1\right)/(2q)\right\} \right]}{2q\exp\left\{ \i\pi(2q-1)\left(2j+1\right)/(2q)\right\} },\;x,t,\in\mathbb{R}.
\label{eq:SSkernel}
\eeq
Here and subsequently $\i=(-1)^{1/2}$. At the second step the boundary conditions are matched, which leads to
the corresponding boundary kernel ${\cal{W}}^b_{ss}(x,t)$. 
An explicit expression for this boundary kernel for $q=1,2$ can be
found, for example in \citet{Agnan96}.
For
$q>2$ derivation of ${\cal{W}}^b_{ss}(x,t)$ becomes very tedious. 

For
regression splines the asymptotic equivalent kernel on $\mathbb{R}$ ${\cal{W}}_{rs}(x,t)$ was obtained by
\citet{HuangStudden93} for $p=3$ only as a $L_2$ projection kernel on a
spline space defined on $\mathbb{R}$. More precisely, for $f\in L_2(\mathbb{R})$,
$$
\int_{-\infty}^\infty{\cal{W}}_{rs}(\cdot,t)f(t)dt=\arg\min_{s\in{{\cal{S}}(p,\mathbb{Z})}}\int_{-\infty}^\infty\{s(x)-f(x)\}^2dx.
$$
The spline space ${\cal{S}}(p,\mathbb{Z})=\{s(x):\;s(x)=\sum_{i=-\infty}^\infty
N_i(x)\theta_i,\;\theta=(\theta_i)_{i\in\mathbb{Z}}\in\ell_2\}$, where $N_i(x)$ denotes a B-spline centered at $i$,  is the
$L_2(\mathbb{R})$ subspace of splines with integer knots. 
With $h=K^{-1}$ \citet{HuangStudden93} give the expression for
$h{\cal{W}}_{rs}\left(hx,ht\right)={\cal{K}}_{rs}(x,t)$ in terms of
normalised cubic B-splines. In contrast to the asymptotic equivalent
kernel for smoothing splines, ${\cal{K}}_{rs}(x,t)$ is not translation-invariant. The boundary kernel for regression splines ${\cal{W}}^b_{rs}(x,t)$ has
not been obtained.

The derivations of ${\cal{W}}_{ss}(x,t)$ and ${\cal{W}}_{rs}(x,t)$ seem to differ greatly both technically and
conceptually. In the following we derive the
asymptotic equivalent kernel on the real line ${\cal{W}}(x,t)$ for
all spline estimators in a unified framework for general $p$, $q$,
$K$ and $\lambda$, and study the pointwise asymptotic
properties of all spline estimators away from the boundaries.
\section{Periodic spline spaces}
\label{seq:periodic}

Let us first introduce some notation. Let
\beq
Q_{p-1}(z)=\sum_{l=-\infty}^\infty\sinc\{\pi(z+l)\}^{p+1},
\label{eq:Q_polynom}
\eeq
a polynomial of $\cos(\pi z)$ of degree
$(p-1)$, which can be expressed in terms of Euler--Frobenius
polynomials $\Pi_p(\cdot)$ \citep[][]{Schoenberg73}: 
\beqn
Q_{p-1}(z)&=&\begin{cases}
\exp\{\i z\pi(p-1)\}{\Pi}_p\{\exp(-2\i\pi z)\}/p!,&p\mbox{ odd},\\
\exp\{\i z\pi(p-1)\}{\widetilde{\Pi}}_p(z)/p!,&p\mbox{ even},\end{cases}
\eeqn
for 
$$\widetilde{\Pi}_p(z)=\cos\left(\pi{z}/{2}\right)^{p+1}\Pi_{p}\left\{\exp(-\pi\i
      z)\right\}-(-1)^{p/2}\i\sin\left(\pi{z}/{2}\right)^{p+1}\Pi_{p}\left\{-\exp(-\pi
      iz)\right\}.$$
Further, we make use of the exponential splines
\citep[][]{Schoenberg73}
\beq
\Phi_p(t,z)={z^{\left\lfloor t\right\rfloor
  }}\left(1-z^{-1}\right)^p\sum_{j=0}^p\binom{p}{j}\frac{{\left\{
      t\right\}}^{p-j}\:
  {\Pi}_j\left(z\right)}{p!\left(z-1\right)^j},\;\;z\neq 0,\;z\neq 1,
\label{eq:Phi_splines}
\eeq
where $\{t\}$ denotes the fractional part of $t$ and
$\lfloor{t}\rfloor$ is the largest integer not greater than $t$. With
the convention $0^0=1$, one can also define $\Phi_p(t,1)=1$. Note that
\beq 
\Phi_{p}\{t ,\exp(2\pi \i z)\}=\frac{\exp(2\pi\i zt)}{\exp\{\pi\i z(p+1)\}}\sum_{l=-\infty}^{\infty}(-1)^{l(p+1)}\sinc\left\{
  \pi\left(z+l\right)\right\}^{p+1}\exp(2\pi\i lt).
\label{eq:F_polynom} 
\eeq 
Next, we define
\beq
Q_{p,\scriptsize{M}}\left(z\right)=\frac{1}{N}\sum_{i=1}^{N}\left|\Phi_{p}\left\{i/M+(p+1)/2,\exp(-2\pi\i z)\right\}\right|^2,
\label{eq:Qpm}
\eeq         
where $M=N/K$. For $M=1$ we find $Q_{p,1}(z)=Q_{p-1}^2(z)$, since $\Phi\left\{ \left(p+1\right)/2 ,\exp\left(2\pi\i
    z\right)\right\}=Q_{p-1}\left(z\right)$  from
(\ref{eq:F_polynom}).
 If
$M=N/K>1$, then $Q_{p,\scriptsize{M}}(z)$ varies between $Q_{2p}(z)$
and $Q^2_{p-1}(z)$, depending on $M$. In particular, it can be shown
that $Q_{p,\scriptsize{M}}(z)=Q_{2p}(z)+c\sin(\pi
z)^{p+1}M^{-(p+1)}$, for a constant $c>0$.

Now we state a lemma giving the explicit expression for
the complex-valued \citet{DemmlerReinsch75} basis for the periodic
spline space ${\mathcal{S}}_{{\tt
    per}}(p;\sett)=\{s:s\in{\mathcal{S}}(p;\sett)\mbox{  and  }s^{(j)}(0)=s^{(j)}(1),\;j=0,\ldots,p-1\}$.
\begin{lemma}
\label{lemma1}
The functions
\beq
\psi_{i}\left(x\right)=\frac{\Phi_p\{Kx+(p+1)/2,\exp(-2\pi\i
i/K)\}}{\left\{Q_{p,M}(i/K)\right\}^{1/2}}\;\; (i=1,\ldots,K), \;\;x\in\mathbb{R},
\label{eq:basis_b}
\eeq
form the complex-valued Demmler-Reinsch basis in
${\mathcal{S}}_{{\tt per}}(p;\sett)$, i.e.,
\beq
\frac{1}{N}\sum_{i=1}^N\psi_i(l/N)\overline{\psi_j(l/N)}&=&\delta_{i,j},
\label{eq:DRcond1}\\
\int_0^1\psi^{(q)}_i(x)\overline{\psi^{(q)}_j(x)}dx&=&\nu_i\delta_{i,j}\;\;(i,j=1,\ldots, K)
\label{eq:DRcond2}
\eeq
and the eigenvalues
\beq
\nu_i=(2\pi i)^{2q}\mathrm{sinc}(\pi i/K)^{2q}\frac{Q_{2p-2q}(i/K)}{Q_{p,M}(i/K)}.
\label{eq:eigenvalues}
\eeq
Moreover, the functions
\beqn
\phi_{i}\left(x\right)=\frac{\Phi_p\{Kx+(p+1)/2,\exp(-2\pi\i
i/K)\}}{\left\{Q_{2p}(i/K)\right\}^{1/2}} \;\;(i=1,\ldots,K),
\eeqn
satisfy 
\beq
\int_0^1\phi_i(x)\overline{\phi_j(x)}dx=\delta_{i,j}=
\frac{1}{\mu_i}\int_0^1\phi^{(q)}_i(x)\overline{\phi^{(q)}_j(x)}dx
\label{eq:DRc}
\eeq
for $\mu_i=(2\pi i)^{2q}\mathrm{sinc}(\pi i/K)^{2q}{Q_{2p-2q}(i/K)}/{Q_{2p}(i/K)}$.
\end{lemma}

The basis functions $\psi_i(x)$ is the scaled
  discrete Fourier transform of periodic B-splines 
    $ \psi_i(x)K\{Q_{p,M}(i/K)\}^{1/2}=\sum_{l=1}^KB_l(x)\exp(-2\pi\i il/K)$.
    A similar basis up to a scaling factor for $N=K$ has been considered by
    \citet{Lee-etal92} and \citet{Zheludev98}.
    
Even though the Demmler--Reinsch basis for periodic smoothing splines
  was employed by \citet{CogburnDavis74} and
  \citet{CravenWahba78}, no explicit expressions for $\psi_i$ and
  $\nu_i$ were given there. For $K=N$ and $p=2q-1$, 
$\nu_i=(2\pi\i)^{2q}\sinc(\pi i/K)^{2q}Q_{2q-2}(i/K)^{-1}$ and at the data
points $l/N$, the Demmler--Reinsh basis reduces to $\psi_i(l/N)=\exp(-2\pi\i il)$.

Thus, any $s\in {\mathcal{S}}_{{\tt per}}(p;\sett)$ can be
represented as $s(x)=\sum_{i=1}^{K}\beta_i\psi_i(x)$ and the solution
to (\ref{eq:min_problem}) over the special class of periodic splines
${\mathcal{S}}_{{\tt per}}(p;\sett)$ results in 
$$\widehat{f}_{\tt
  per}(x)=\frac{1}{N}\sum_{i=1}^N\sum_{j=1}^{K}\frac{\psi_j(x)\overline{\psi_j(x_i)}}{1+\lambda\nu_j}Y_i=\frac{1}{N}\sum_{i=1}^NW_{\tt
  per}^{[0,1]}(x,x_i)Y_i,$$ 
where 
$
{W}_{\tt per}^{[0,1]}(x,t)=\sum_{i=1}^K\psi_i(x) {\overline{\psi_i(t)}}{(1+\lambda\nu_i)^{-1}}
$, with ${\overline{\psi_i(t)}}$ denoting the complex conjugate, 
is the effective kernel for periodic spline estimators, which depends
on $N$ via $Q_{p,M}$. The corresponding asymptotic equivalent kernel is 
$
{\cal{W}}_{\tt per}^{[0,1]}(x,t)=\sum_{i=1}^{K}{\phi_i(x)\overline{\phi_i(t)}}{(1+\lambda\mu_i)^{-1}}
$,
and is such that
\beqn
\int_0^1{\cal{W}}_{\tt
  per}^{[0,1]}(\cdot,x)f(x)dx=\arg\min_{s\in{\mathcal{S}}_{{\tt
      per}}(p;\sett)}\left[\int_0^1\{s(x)-f(x)\}^2dx+\lambda\int_0^1\{s^{(q)}(x)\}^2dx  \right].
\eeqn
Both ${W}_{\tt per}^{[0,1]}(x,t)$ and
${\cal{W}}_{\tt per}^{[0,1]}(x,t)$ are known explicitly.


\section{Asymptotic equivalent kernels on $\mathbb{R}$}
\label{sec:EK}
Let
$$
{\mathcal{W}}(x,t)=\int_0^K \frac{\phi(u,x)\overline{\phi(u,t)}}{1+\lambda\mu(u)}du,
$$
where $\mu(u)=(2\pi u)^{2q}\sinc(\pi
u/K)^{2q}Q_{2p-2q}(u/K)/Q_{2p}(u/K)$ and
$\phi(u,x)=\Phi_p\{Kx+(p+1)/2,\exp(-2\pi \i
u/K)\}\{Q_{2p}(u/K)\}^{-1/2}$. 
As shown in the proof of Lemma \ref{lemma2} ${\cal{W}}_{\tt per}^{[0,1]}(x,t)$ can be obtained by folding back $
{\mathcal{W}}(x,t)$, that is
$
{\cal{W}}_{\tt per}^{[0,1]}(x,t)=\sum_{l=-\infty}^\infty
{\mathcal{W}}(x,t+l)
$. In particular, for a periodic function $f$ one finds
$
\int_0^1{\cal{W}}_{\tt per}^{[0,1]}(x,t)f(t)dt=\int_{-\infty}^\infty{\cal {W}}(x,t)f(t)dt
$. 
Subsequently, we refer to ${\cal{W}}(x,t)$ as the {asymptotic
  equivalent kernel for spline estimators on $\mathbb{R}$}. The following lemma gives the explicit expression for ${\cal{W}}(x,t)$.
\begin{lemma}
\label{lemma2}
Let 
\beq
P_{2p}(u)={\Pi}_{2p}(u)+(-1)^q\lambda
K^{2q}(1-u)^{2q}\Pi_{2p-2q+1}(u)/(2p-2q+1)!
\label{eq:P2p}
\eeq
be a polynomial of degree $2p$ where $ \Pi_{p}(u)$ is the
Euler--Frobenius polynomial. Let also $r_j$, $r_j^{-1}$,
$j=1,\ldots,p$ be the roots of  $P_{2p}(u)$ with $|r_j|<1$. Then, denoting $P^{'}_{2p}(r_j)=\left.\partial
  P_{2p}(u)/\partial u\right|_{u=r_j}$, $d_{x,t}=\left\lfloor Kx-\left\{ (p+1)/2\right\} \right\rfloor -\left\lfloor Kt-\left\{ (p+1)/2\right\} \right\rfloor  $ and representing
$
z^{p-d_{x,t}}\Phi_p\{Kx+(p+1)/2,z
\}\Phi_p\{Kt+(p+1)/2,z^{-1}\}=\sum_{l=0}^{2p}\alpha_l(\{Kx\},\{Kt\})z^{l}
$
for some functions $\alpha_l(t_1,t_2)$ and $x,t\in\mathbb{R}$, results in
\beq
\mathcal{W}(x,t)=K\sum_{j=1}^p\sum_{l=0}^{2p}\frac{\alpha_l\left(\{Kx\},\{Kt\}
  \right)}{P_{2p}^{'}(r_j)}r_{j}^{\left|d_{x,t}+l-1\right|+(2p-2)\mathbb{I}(d_{x,t}+l\leq0)},
\label{eq:kernelW}
\eeq
where $\mathbb{I}(A)$ is an indicator function, which is
equal to $1$ if $A$ is true and $0$ otherwise.
\end{lemma}

 For $p=q=1$ equivalent kernel ${\cal{W}}(x,t)$ has a simple
representation, for larger $p$ and $q$ it becomes much more involved. For $p=q=1$ one finds $r_1=1-\{(3+36\lambda
  K^2)^{1/2}-3\}/(6\lambda K^2-1)$, $P_{2}^{'}(r_1)=\{(1+12\lambda
  K^2)/3\}^{1/2}$, $\alpha_2(t_1,t_2)=\alpha_0(t_2,t_1)=t_1-t_1t_2$ and
$\alpha_1(t_1,t_2)=1-\alpha_0(t_1,t_2)-\alpha_2(t_1,t_2)$. Expression (\ref{eq:kernelW}) for the asymptotic equivalent kernel on $\mathbb{R}$ is valid for any combinations of
$p$, $q$, $K$ and $\lambda$. The next lemma gives the order of this
kernel. Let
$\mu_{m}({\cal{W}})=\int_{-\infty}^\infty t^m{\cal{W}}(x,t)dt$. 
\begin{lemma}
For $x,t\in\mathbb{R}$, $d=\min\{p+1,2q\}$ and
$m\in\mathbb{N}_0$,
\beqn
\mu_m(\cal{W})&=&\begin{cases}
1,&m=0,\\
0,&m=1,\ldots,d-1,\\
-\mathbb{I}(
    d=p+1)\frac{\mathcal{B}_{p+1}\left(\left\{
        Kx+\frac{p+1}{2}\right\} \right)}{K^{p+1}} -\mathbb{I}(d=2q){\left(-1\right)^{q}}{\lambda\left(2q\right)!},&m=d,\end{cases}
\eeqn
with $\mathcal{B}_{p+1}(x)$ a $(p+1)$-th degree Bernoulli
polynomial.
\label{lemma4}
\end{lemma}
To establish the correspondence between our results and the known asymptotic
equivalent kernels on $\mathbb{R}$ for smoothing and regression
spline estimators, we must first determine an appropriate bandwidth,
which is universal for all spline estimators. 
Let us define the variable $k_q=\lambda^{1/(2q)}\pi K$, which characterizes the type of the spline
estimator. In particular, $k_q=0$ corresponds to the regression spline
estimator, $k_q=\lambda^{1/(2q)}\pi N$ to the smoothing
spline estimator and all intermediate values characterize penalized
spline estimators. 
With this, we introduce the bandwidth $h(k_q)$, which is universal
for all spline estimators, given by 
$$
h(k_q)^{-1}=\int_0^{K}\frac{dx}{1+\lambda (\pi x)^{2q}}=\lambda^{-1/(2q)}\pi^{-1} \int_0^{k_q}\frac{dx}{1+x^{2q}}.
$$
Bandwidth $h(k_q)$ is a smooth function of $k_q$ with a rather
complicated closed form
expression available for each $q$. In our subsequent
developments we use the following representation. 
\beq
h(k_q)^{-1}=\lambda^{-1/(2q)}\pi^{-1}\begin{cases}
k_q\:c_1,&k_q<1,\\
\pi\: c_2,&k_q\geq1,
\end{cases}\label{eq:bandwidth}
\eeq
with constants
$c_{1}=c_1(k_q)={\leftsup{2}{{\cal{F}}}_1[\{1,1/(2q)\};\{1+1/(2q)\},-k_q^{2q}]}$
and $ c_{2} = c_2(k_q)=\widetilde{c}_2
-\pi^{-1}{k_{q}}^{1-2q}\leftsup{2}{\cal{F}}_1[\{1,1-1/(2q)\};\{2-1/(2q)\},-k_q^{-2q}]/(2q-1)
$, where $\widetilde{c}_2=\pi^{-1}{\sinc \{ \pi /(2q)\}}^{-1}$is independent of $k_q$ and
$\leftsup{2}{\cal{F}}_1$ denoting the hypergeometric series
\citep[][]{AbramowitzStegun72}. Both $c_{1}(k_{q})$ and
$c_{2}(k_{q})$ 
are convergent and vary slowly with $k_q$, namely
$c_1(k_q)\in(\pi/4,1]$ for any $k_q<1$ and $c_2(k_q)\in(1/4,1/2]$ for any
$k_q\geq 1$. For the regression
spline estimators, $k_q=0$, the bandwidth
$h(0)=K^{-1}$ and for the smoothing spline estimators, $k_q\rightarrow\infty$, the bandwidth
$h(\infty)=\lambda^{1/(2q)}/\widetilde{c}_2$. 

This separation into two cases for $k_q<1$ and $k_q\geq 1$ was introduced by \citet{Claeskens-etal09}. If $k_q<1$, then the asymptotic behaviour of the penalized spline
estimator is similar to that of the regression spline estimator, while
$k_q\geq 1$ corresponds to asymptotic behaviour similar to that of
smoothing splines. 

\begin{figure}[!t]
\begin{center}
\begin{tabular}{cc}
\includegraphics[width=0.5\textwidth]{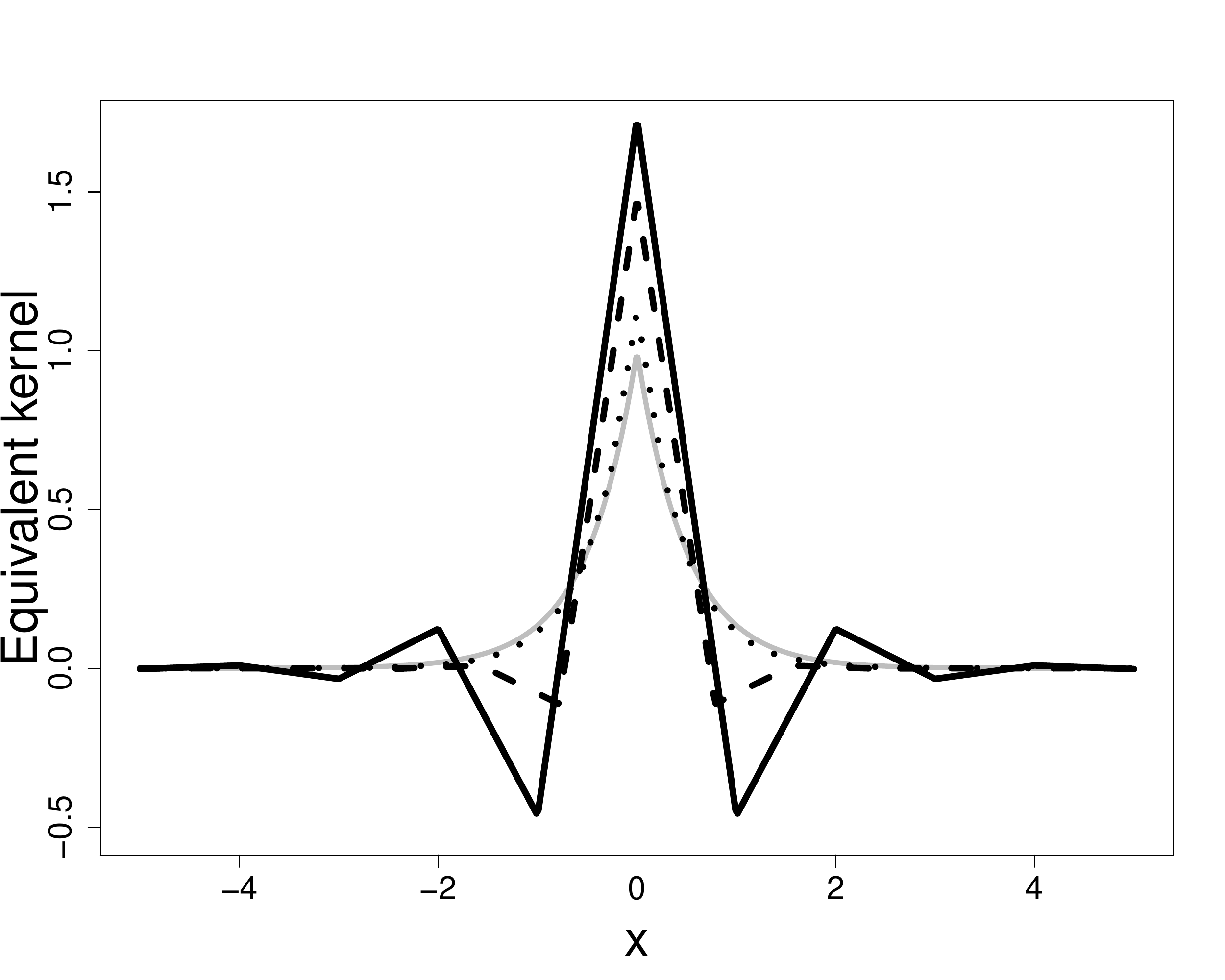}
&\includegraphics[width=0.5\textwidth]{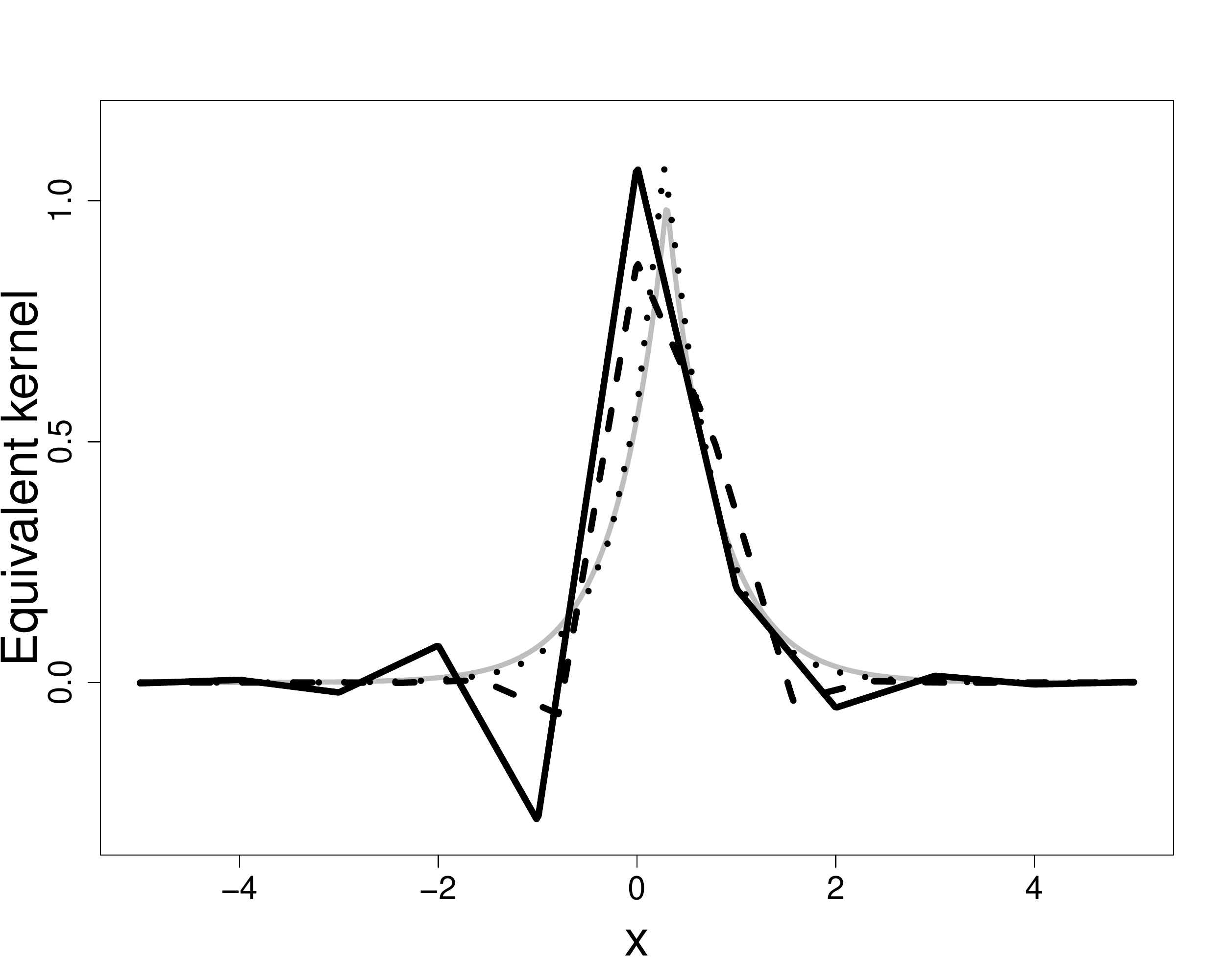}\\
\includegraphics[width=0.5\textwidth]{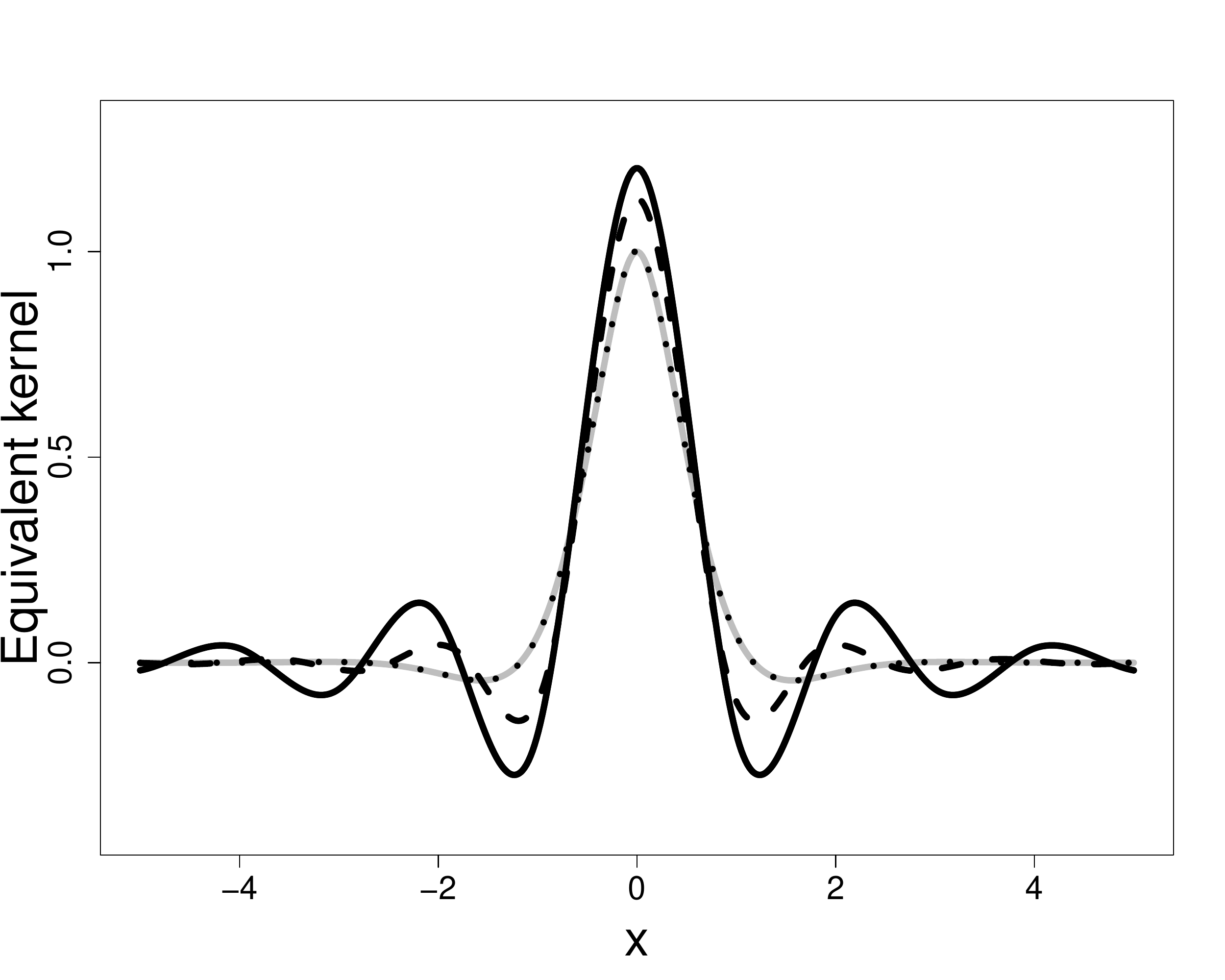}
 &\includegraphics[width=0.5\textwidth]{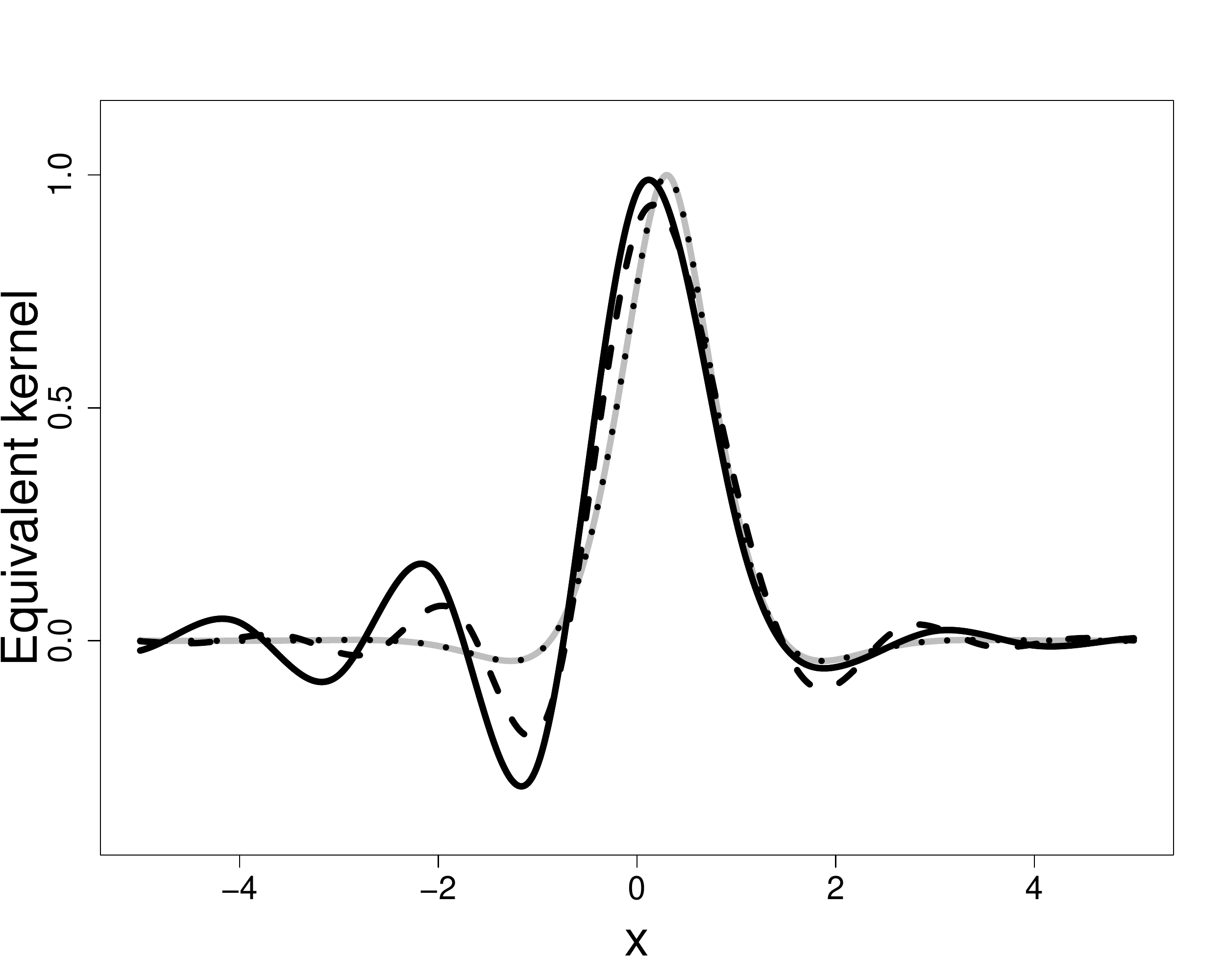}\\\vspace{0.4cm}
\end{tabular}
\caption{Equivalent kernels ${\cal{K}}(x,0)$ (left) and
  ${\cal{K}}(x,0.3)$ (right) for $p=q=1$ (top) and
  $p=2q-1=3$ (bottom) and different values of $k_q$. The grey line corresponds to the smoothing
  spline kernel, the bold line to the regression spline kernel, the dashed and dotted lines to the penalized spline kernel with $k_q=1$ and $k_q=5$, respectively.}
\label{fig:psKernel}
\end{center}
\end{figure}

Let $\mathcal{K}(x,t)$  denote the scaled version of the
equivalent kernel on
$\mathbb{R}$, i.e.,
$$h(k_q)\mathcal{W}\left\{h(k_q)x,h(k_q)t\right\}=\mathcal{K}(x,t).$$
\begin{theorem}
\label{theorem2}
The equivalent kernel for spline estimators on
$\mathbb{R}$ with $p=2q-1$ satisfies
$$
\begin{cases}
c_{1}\mathcal{K}\left(c_{1}x,c_{1}t\right)=\mathcal{K}_{rs}(x,t)-k_{q}^{2q}\mathcal{K}_{1}(x,t)
, & k_{q}<1,\\
c_{2}\mathcal{K}\left(c_{2}x,c_{2}t\right)=\mathcal{K}_{ss}(t-x)+k_{q}^{1-2q}\mathcal{K}_{2}(x,t), & k_{q}\geq1,
\end{cases}
$$ 
where  $\mathcal{K}_{1}(x,t)$ and $\mathcal{K}_{2}(x,t)$
are bounded functions 
given in the proof, $\mathcal{K}_{{rs}}(x,t)$ is the asymptotic regression spline
equivalent kernel on $\mathbb{R}$, i.e.,
\beqn
\mathcal{K}_{{rs}}(x,t)=\sum_{j=1}^p\sum_{l=0}^{2p}\frac{\alpha_l\left(\{x\},\{t\}
  \right)}{P_{2p}^{'}(r_j)}\;r_{j}^{\left|d_{x,t}+l-1\right|+(2p-2)\mathbb{I}(d_{x,t}+l\leq0)},
\eeqn
with $\alpha_l$, $r_j$, $d_{x,t}$, $P_{2p}(u)=\Pi_{2p}(u)$ defined in
Lemma \ref{lemma2} and $\mathcal{K}_{{ss}}(x,t)$ the asymptotic
smoothing spline equivalent kernel on $\mathbb{R}$ as given in (\ref{eq:SSkernel}).
\end{theorem}
Theorem \ref{theorem2} implies that if $k_q<1$, then the
asymptotic equivalent kernel on $\mathbb{R}$ for spline estimators is dominated by
${\cal{K}}_{rs}(x,t)$, the asymptotic equivalent
regression spline kernel on $\mathbb{R}$, while for $k_q\geq 1$ ${\cal{K}}(x,t)$ is dominated
by ${\cal{K}}_{ss}(x,t)$, which agrees with the findings of
\citet{Claeskens-etal09}. Moreover, if $p=2q-1$, then
${\cal{K}}(x,t)$ varies smoothly
between ${\cal{K}}_{rs}(x,t)$ and ${\cal{K}}_{ss}(x,t)$, appropriately
scaled, both having the same order. In particular, $\lim_{k_q\rightarrow\infty}{c_2\cal{K}}(c_2x,c_2t)=\widetilde{c}_2^{-1}{\cal{K}}_{ss}\left\{(x-t)/\widetilde{c}_2\right\}$
and $\lim_{k_q\rightarrow
  0}c_1\mathcal{K}(c_1x,c_1t)=
\mathcal{K}_{rs}(x,t)$. 
Figure
\ref{fig:psKernel} depicts the penalized spline kernel
${\cal{K}}(x,t)$ at $t=0$ and $t=0.3$ as a function of $x$ for different values
of $k_q$ and for $p=1,3$. The case $k_q=0$ corresponds to
${\cal{K}}_{rs}(x,t)$. As $k_q$ grows, ${\cal{K}}(x,t)$ becomes
more symmetric and for $k_q=5$ is already
indistinguishable from the smoothing spline kernel.

From Lemma \ref{lemma2} the expression for ${\cal{K}}(x,t)$ can
be obtained for any combination of $p$ and $q$. However, if $p\neq
2q-1$, then no smooth transition between two scenarios with $k_q<1$ and
$k_q\geq 1$ is possible. In this case the orders of the asymptotic
equivalent regression and
smoothing spline kernels will be different and, hence, the asymptotic
rates of the corresponding estimators, see also discussion on the
parameter choice in Section \ref{sec:SmoothPar}. To keep the exposition clear,
in the next section we focus on the case $p=2q-1$.


\section{Local asymptotics of spline estimators}
\label{sec:local}
The following theorem gives the pointwise bias and variance of all
spline estimators away from the boundaries.
\begin{theorem}
\label{theorem3}
Let the model (\ref{eq:model}) hold and $\widehat{f}(x)\in
{\mathcal{S}}(2q-1;\sett)$ be the solution to (\ref{eq:min_problem})
with $x_i=i/N$ ($i=1,\ldots, N$) and
$\underline{\tau}_{\ss{K}}=\{i/K\}_{i=0}^K$. Then for  $f\in
\mathcal{H}_{2q}[0,1]$, such that $f^{(2q)}$ is
H\"older continuous with $|f^{(2q)}(x)-f^{(2q)}(t)|\leq
L|x-t|^\alpha$, $x,t\in [0,1]$, $L>0$ and $\alpha\in(0,1]$, 
for any $x\in{\mathcal{I}}_q=\left[\delta_q h(k_q)\log\left\{h(k_q)^{-1}\right\},1-\delta_q h(k_q)\log\left\{h(k_q)^{-1}\right\}\right]$, we have
\beq
{E}\left\{ \widehat{f}(x)\right\} -f(x) & =&(-1)^{q+1}h(k_q)^{2q}\frac{f^{(2q)}(x)}{(2q)!}C(k_{q},x)+o\left\{ h(k_{q})^{2q}\right\}\label{eq:bias}
 \\
\mathrm{var}\left\{\widehat{f}(x)\right\} & =&\frac{\sigma^{2}}{Nh(k_{q})}\int_{-\infty}^{\infty}{\cal K}\left\{x/h(k_q),t\right\}^2dt+o\left\{N^{-1}h(k_q)^{-1}\right\},\label{eq:var}
\eeq
where
\beqn
C(k_{q},x)=\begin{cases}
c_{1}^{2q}\left[(-1)^{q}\mathcal{B}_{2q}\left(\left\{ Kx\right\}
  \right)+(2q)!\pi^{-2q}
  k_{q}^{2q}\right],
& k_q<1\\
c_2^{2q}\left[(2q)! +(-1)^{q}\mathcal{B}_{2q}\left(\left\{ Kx\right\} \right) \pi^{2q} k_{q}^{-2q}\right],&k_q\geq 1\end{cases}
\eeqn
and $\int_{-\infty}^{\infty}{\cal
  K}\left\{x/h(k_q),t\right\}^2dt<C^2/\log(\gamma^{-1})$, for some
$C\in(0,\infty)$ and $\gamma\in(0,1)$, both depending on $k_q$, explicitly given in the proof of
Lemma 4 in the Appendix. The constant $\delta_q$ in
${\cal{I}}_q$ depends only on $q$ and $\gamma$
and it is such that $\delta_q>2q\log_\gamma(1/e)$.

If $f$ is a periodic function, that is $f\in
{\mathcal{P}}_{2q}[0,1]=\{f:\;f\in
C^{2q}(\mathbb{R}),\;f^{(j)}(0+l)=f^{(j)}(1+l),\;l\in\mathbb{Z},\;j=0,\ldots,2q-1\}$
and $\widehat{f}(x)=\widehat{f}_{\tt per}(x)\in
{\mathcal{S}}_{\tt per}(2q-1;\sett)$, then (\ref{eq:bias}) and (\ref{eq:var}) hold uniformly
for all $x\in[0,1]$.
\end{theorem}

It is easy to see that the limiting cases of expression (\ref{eq:bias}) and (\ref{eq:var}) for $k_q=0$ and $k_q\rightarrow\infty$ coincide with the known results on regression and smoothing splines. 
\cite{LiRuppert07} obtained pointwise bias and variance of penalized spline estimators with a slightly different penalty matrix for the special cases $p=0,1$ with $q=1,2$ and $k_q>1$. These results also agree with the corresponding equations (\ref{eq:bias}) and (\ref{eq:var}), up to the expression for the bandwidth, which is given up to a constant only.

An asymptotic optimal bandwidth at any
$x\in{\cal{I}}_q$ can be obtained from Theorem \ref{theorem3}. 
\begin{corollary} \label{lemma:improvedAsymptotic}
Under the assumptions of Theorem \ref{theorem3}
the asymptotic optimal bandwidth depending on $k_q$ for
$f\in{\cal{H}}_{2q}[0,1]$ at any
$x\in{\cal{I}}_q$ and for $f\in{\cal{P}}_{2q}[0,1]$ at any $x\in[0,1]$ is 
$$h_{opt}(k_q,x)=\left[ N \frac{4qC(k_q,x)^2\left\{f^{(2q)}(x) \right\}^2}{\sigma^2(2q!)^2\int_{-\infty}^{\infty}{\cal K}\left\{x/h(k_q),t\right\}^2dt}\right]^{-1/(4q+1)}.$$ 
\end{corollary}
The proof is straightforward from (\ref{eq:bias}) and (\ref{eq:var}).

\section{Choice of parameters}
\label{sec:SmoothPar}
Several parameters for penalized spline
estimators must be
chosen in practice, i.e., $p$, $q$, $K$ and $\lambda$.
Theorems \ref{theorem2} and \ref{theorem3} allow us 
to make the following practical recommendations. 

Theorem \ref{theorem2} implies that setting $p=2q-1$
provides a smooth transition between two asymptotic scenarios, that is $k_q<1$
and $k_q\geq 1$. In this case, according to Theorem \ref{theorem3},
the convergence rate of spline estimators in both scenarios is the same and $k_q$ enters
only the constants, see definition of $C(k_q,x)$ in Theorem
\ref{theorem3}. Hence, it is convenient to choose $q$ and set $p=2q-1$ to make
the rate of convergence of the estimator independent on the number of
knots chosen. The choice of $q$ depends on the
smoothness of the underlying regression function $f$, but $q=2$ is typically
taken in practice. 

A second issue is the choice of tuning parameters. 
For regression or smoothing splines, there is only one tuning
parameter, $K$ or $\lambda$, respectively, which is typically chosen to minimize an unbiased estimator
of the empirical $L_2$ risk of the estimator, for example the generalized cross
validation criterion. However, in practice one
typically prefers to use penalized splines with $K\ll N$, but still with a penalty tuned by $\lambda$, so that two parameters need to be selected. A typical approach in practice is to fix $K$ arbitrarly and then choose $\lambda$. \citet{Ruppert02} ran a
large simulation study, recommending using $\min\{N/4,35\}$
knots in practice. The results of Theorem \ref{theorem3} suggest that first
fixing $K$  and then choosing $\lambda$ can be problematic. Indeed, if $K$ is selected so that
$k_q<1$, then the bandwidth $h(k_q)$ would depend on
$K$ with $\lambda$ entering only $c_1(k_q)$. Hence, in this case
$\lambda$ cannot be estimated consistently. We argue that both
parameters should be chosen simultaneously. One can fix a reasonable
value of $k_q$ and search only over those $K$s and $\lambda$s that give
this particular $k_q$. Since for $p=2q-1$ the convergence rate is
independent of $k_q$, the difference
between estimators with different $k_q$ values would vanish with the growing sample size. The constant
$C(k_q,x)$ from Theorem \ref{theorem3} is a smooth increasing function
of $k_q$, so smaller values could be
preferable. However, $k_q<1$ values lead to the regression spline-type
estimators, for which not only the number of knots but also their locations
is crucial. Hence, fixing $k_q$ slightly larger than unity in practice would
help to avoid the dependence on knot location.

To illustrate this discussion we ran a simulation study with two functions, $f_1(x)=\sin(6\pi x)$ and $f_2(x)=\sin(2\pi x)^2\exp(x)$. We
set $n\in\{300,1000\}$, $\sigma=0.1$, $p=3$, $q=2$ and
consider $k_q\in\{0.5,1,1.2,1.5,5\}$. We restricted the range of $K$
to $K=2,\ldots,50$ and set $\lambda=\{k_q/(\pi
K)\}^{2q}$ for each given $k_q$ and all values $K$. Finally, we evaluated the generalized cross validation criterion $\mbox{GCV}(k_q)$, $k_q\in\{0.5,1,1.2,1.5,5\}$,
for each $k_q$ choosing those $K$ and $\lambda$, that minimize
$\mbox{GCV}(k_q)$. This procedure is very fast, since the only values of $K$ and $\lambda$ 
considered are such that $\lambda^{1/(2q)}\pi K$ equals a particular number, which
leads to a search over a very sparse grid. Table
\ref{Table1} reports the results from $M=500$ Monte Carlo simulations;
here
$A_N(f)=(NM)^{-1}\sum_{j=1}^N\sum_{i=1}^{M}\{\widehat{f}_i(x_j)-f(x_j)\}^2$,
with $\widehat{f}_i(\cdot)$ denoting the estimator of $f$ in the $i$-th Monte Carlo
replication. 
{\renewcommand{\arraystretch}{1.2} 
\begin{table}
\label{Table1}
\caption{$A_N(f_1)$ and $A_N(f_2)$
  depending on $N$ and $k_q$. All entries are multiplied by $10^4$.} 
\begin{center}
\begin{tabular}{ccccccc}
   &   &  $k_q=0.5$ & $k_q=1$ & $k_q=1.2$ &$k_q=1.5$ & $k_q=5$ \\
$N=300$&$A_N(f_1) $&6.781 &6.720&  6.662& { 6.519}&6.798\\
               &$A_N(f_2) $& 5.212 &  5.158& {  5.139}& 5.206&5.896\\
$N=1000$&$A_N(f_1) $&2.124 & 2.119 &{2.113}  &{2.098} &2.253\\
               &$A_N(f_2) $&1.564  &1.560  &{ 1.556}  &1.562 &2.018
\end{tabular}
\end{center}
\end{table}
Choosing $k_q$ slightly larger than unity leads to a somewhat better average
mean squared error, but the differences
between spline estimators with different $k_q$s are marginal.

\section{Discussion}
\label{sec:extentions}
We have obtained equivalent kernels and local asymptotic results for
spline estimators of sufficiently smooth $f$s away from the
boundaries, assuming equidistant knots and
observations. For periodic functions all the results hold uniformly
over $[0,1]$. Two issues remain undiscussed. First, non-equidistant design for knots and
observations and second, the boundary behaviour of spline estimators of
general smooth functions. 

The equidistant design assumption is dominant in the literature on equivalent
kernels, since it allows a clear exposition. However, it can easily be relaxed, as in \citet{HuangStudden93}. In particular, if the design points $x_i$ ($i=1,\ldots,N$)
have a limiting density $g(x)$ and the sequence of knots
$\underline\tau_K$ is such that $\int_{\tau_{i-1}}^{\tau_i}p(x)dx=1/K$, for a positive
continuous density $p(x)$ on $[0,1]$, then the equivalent kernel for a 
general spline estimator satisfies
$$
{\cal{W}}(x,t)=\frac{1}{g(x)h\{k_q(x)\}}{\cal{K}}\left[\frac{x}{h\{k_q(x)\}},\frac{t}{h\{k_q(x)\}}\right],
$$
where $k_q(x)=k_qp(x)^{1-1/(2q)}$ and $h\{k_q(x)\}^{-1}=\int_0^K\left\{1+k_q(x)^{2q}(t/K)^{2q}\right\}^{-1}dt$.
The asymptotic results of Theorem \ref{theorem3} should then be read
as follows: $k_q$ is everywhere replaced by $k_q(x)$ and the variance in (\ref{eq:var}) will be additionally
scaled by $1/g(x)$.

The boundary behaviour of general
spline estimators remains open. It can be studied after derivation of a boundary kernel
${\cal{W}}^b(x,t)$, such that the equivalent kernel on $[0,1]$ is ${\cal{W}}^{[0,1]}(x,t)={\cal{W}}(x,t)+{\cal{W}}^b(x,t)$. While such a boundary kernel is
known for smoothing splines, it is not available for regression and penalized spline
estimators, since the Green's function approach of smoothing splines
cannot be applied. In general, it is known that regression spline
estimators do not have boundary effects, while smoothing spline estimators
have a larger bias at the boundaries. Considering the results of
Theorem \ref{theorem2} one can make a
conjecture on the boundary behaviour of general spline estimators. Since ${\cal{K}}(x,t)$ varies smoothly between
${\cal{K}}_{rs}(x,t)$ and ${\cal{K}}_{ss}(x,t)$, one can expect that
additional boundary terms in ${\cal{K}}^{[0,1]}$ also vary
smoothly between ${\cal{K}}_{rs}^b(x,t)$ and ${\cal{K}}_{ss}^b(x,t)$, so
that the boundary effects of spline estimators grow as
$k_q\rightarrow\infty$. This is another reason to select a smaller
$k_q>1$ in practice. 

\begin{appendix}
\label{appendix}

\section*{Appendix. Technical details}
\label{appendix1}

\subsection{Proof of Lemma 1}
Inserting the Fourier series of a periodic B-spline  into the discrete Fourier transform of B-splines, we find
\beqn
\sum_{i=1}^KB_i(x)\exp(-2\pi\i
li/K)&=&\sum_{m=-\infty}^\infty\exp(-2\pi \i mx)\sinc(\pi m/K)^{p+1}\sum_{i=1}^K\exp\{2\pi\i i(m-l)/K\}\\
&=&K\sum_{n=-\infty}^\infty\exp\{-2\pi \i(l+nK)x\}\sinc\{\pi
(l/K+n)\}^{p+1}\\
&=&\{{Q_{p,\scriptsize{M}}(i/K)}\}^{1/2}\psi_l(x),
\eeqn
where in the last equality the representation 
\beq
\psi_i(x)=\{{Q_{p,\scriptsize{M}}(i/K)}\}^{-1/2}\sum_{l=-\infty}^{\infty}\sinc\{\pi(i/K+l)\}^{p+1}\exp\{-2\pi
\i x(i+lK)\}.
\label{eq:DR_sinc}
\eeq
has
been used, which follows from (\ref{eq:F_polynom}), and $n=(m-l)/K$. The properties of the discrete Fourier transform ensure that the
functions $\psi_i(x)$ ($i=1,\ldots,K$) are also the basis in
${\mathcal{S}}(p;\sett)$. 
Property (\ref{eq:DRcond1}) follows immediately from the definition of
$Q_{p,M}(z)$. To show property (\ref{eq:DRcond2}) one can use
again the representation in (\ref{eq:DR_sinc}) to find
\beqn
\frac{\{Q_{p,\scriptsize{M}}(i/K)\}^{1/2}\phi^{(q)}_i(x)}{(-2\pi\i i)^{q}\sinc(\pi i/K)^{q}}&=&\sum_{l=-\infty}^\infty(-1)^{lq}\sinc\{\pi (i/K+l)\}^{p+1-q}\exp\{-2\pi \i x(i+lK)\},
\eeqn
which implies the assertion and proves the lemma. Property
(\ref{eq:DRc}) of functions $\phi_i$ follows similarly from the definition of $Q_{2p}(z)$. 
 \hfill$\square$
\subsection{Proof of Lemma 2}
Let us show that
\beq
{W}(x,t)=\sum_{l=-\infty}^\infty
{\mathcal{W}}(x,t+l).
\label{eq:fold}
\eeq
This can be proved by showing that the Fourier coefficients of both
functions coincide. 
The Fourier coefficients of ${{W}}(x,t)$
as a
function of $t$ at a fixed $x$ can be found from
\beqn
{{W}}(x,t)=\sum_{l=-\infty}^\infty\sum_{i=1}^K
\frac{\sinc\{\pi(i/K+l)\}^{p+1}\phi_i(x)}{\{Q_{2p}(i/K)\}^{1/2}(1+\lambda \mu_i)}\exp\{2\pi\i t(i+lK)\}.
\eeqn 
Since
$Q_{2p}(i/K)=Q_{2p}(i/K+l)$, $\mu_i=\mu_{i+lK}$ and
$\phi_{i}(x)=\phi_{i+lK}(x)$, we obtain
\beq
a_{l}(x)=\frac{\sinc\{\pi(l/K)\}^{p+1}\phi_l(x)}{\{Q_{2p}(l/K)\}^{1/2}(1+\lambda \mu_l)},\;l\in\mathbb{Z},
\label{eq:FourierCoef}
\eeq
for ${{W}}(x,t)=\sum_{l=-\infty}^\infty a_l(x)\exp(2\pi\i
lt)$. From the Poisson summation formula 
\beq
\int_{0}^{1}\sum_{j=-\infty}^{\infty}\mathcal{W}(x,t+j)\exp(-2\pi \i tl)dt=\int_{-\infty}^{\infty}\mathcal{W}(x,t)\exp(-2\pi itl)dt
\label{eq:FourierCoefTransf}
\eeq
follows the equality of $l$th Fourier coefficients of
$\sum_{j=-\infty}^{\infty}\mathcal{W}(x,t+j)$ and of the Fourier
transform of
$\mathcal{W}(x,t)$. Applying the Poisson summation formula again 
we obtain
\beq
\mathcal{W}(x,t)&=&\int_0^K \sum_{l=-\infty}^\infty
\frac{\sinc\{\pi(u/K+l)\}^{p+1}\phi(u,x)}{\{Q_{2p}(u/K)\}^{1/2}\{1+\lambda
  \mu(u)\}}\exp\{2\pi\i t(u+lK)\}du\nonumber\\
&=&\int_{-\infty}^\infty
\frac{\sinc\{\pi(u/K)\}^{p+1}\phi(u,x)}{\{Q_{2p}(u/K)\}^{1/2}\{1+\lambda
  \mu(u)\}}\exp\{2\pi\i t u\}du.\label{eq:Wtransform}
\eeq
From (\ref{eq:FourierCoefTransf}), (\ref{eq:Wtransform}) and the
inverse Fourier transform follows the equality of the Fourier
coefficients of $\sum_{j=-\infty}^{\infty}\mathcal{W}(x,t+j)$  and
$a_l(x)$ in (\ref{eq:FourierCoef}), which proves (\ref{eq:fold}).

Next we aim to represent and $\mathcal{W}(x,t)$ as a ratio
of two polynomials of exponential functions. The basis functions
$\phi_i(x)$ and $\phi(u,x)$, as well as $Q$ polynomials,
can be expressed in terms of the Euler--Frobenius polynomials of
exponential functions, as shown in Section
\ref{sec:pensplines}. 
With this
\beqn
\mathcal{W}(x,t)&=&\int_{0}^{K}\frac{\exp\left(-2\pi\i d_{x,t} u/K\right)\sum_{l=0}^{2p}\alpha_{l}(\{Kx\},\{Kt\})\exp\left(-2\pi\i u l/K\right)}{P_{2p}\left\{\exp\left(-2\pi\mathtt{i}u/K\right)\right\}}du.
\eeqn     
The coefficients of the partial fractional decomposition of $1/P_{2p}$
are
$1/P_{2p}^{'}\left(r_{j}\right)$ and $1/P_{2p}^{'}\left(r_{j}^{-1}\right)$ correspondent to the roots $r_j$ and $r_j^{-1}$ for $j=1,\ldots, p$. From the representation of
$P_{2p}$ as a function of $\cos^{2}\left(\pi i/K\right) =\{\exp\left(-2\pi\i i/K\right)+\exp\left(2\pi\i i/K\right)+2\}/4$ follows that $P_{2p}^{'}(r_{i}^{-1})=-r_{i}^{2-2p}P_{2p}^{'}(r_{i}^{-1})$. Then
$$
\mathcal{W}(x,t)=  \sum_{j=1}^{p}\sum_{l=0}^{2p}\frac{\alpha_{l}(\{Kx\},\{Kt\})}{P_{2p}^{'}\left(r_{j}\right)}\:R(j,l),
$$
for 
\begin{eqnarray*}
R(j,l)&=&\int_{0}^{K}\left[\frac{\exp\left\{ -2\pi\i \left(d_{x,t}+l\right)u/K\right\} }{\exp(-2\pi\i u/K)-r_{j}}-\frac{r_j^{2p-2}\exp\left\{ -2\pi\i \left(d_{x,t}+l\right)u/K\right\} }{\exp(-2\pi\i u/K)-r_{j}^{-1}}\right]du.
\end{eqnarray*}
Solution to $R(j,l)$ follows from the Cauchy integral formula, where the contour
integral is taken counter-clockwise
$$
R(j,l)=\begin{cases}
\frac{K}{2\pi\i}\oint_{\left|z\right|=1}\left(\frac{z^{d_{x,t}+l-1}}{z-r_{j}}-\frac{r_{j}^{2p-2}z^{d_{x,t}+l-1}}{z-r_{j}^{-1}}\right)dz=Kr_{j}^{d_{x,t}+l-1}, & \left(d_{x,t}+l\right)>0\\
\frac{K}{2\pi\i}\oint_{\left|z\right|=1}\left(\frac{r_{j}^{2p-1}z^{-d_{x,t}-l}}{z-r_{j}}-\frac{r_{j}^{-1}z^{-d_{x,t}-l}}{z-r_{j}^{-1}}\right)dz=Kz^{-d_{x,t}-l+2p-1}, & \left(d_{x,t}+l\right)\leq 0.
\end{cases}
$$
\hfill$\square$
\subsection{Proof of Lemma \ref{lemma4}}
From (\ref{eq:Wtransform}) and symmetry of the kernel
follows that
$
\mathcal{W}(x,t)=\int_{-\infty}^{\infty}\overline{a(u,x)}\exp\left(-2\pi itu\right)du,
$
with $a(u,x)$ defined as 
\beqn
a\left(u,x\right)&=&\frac{\sinc\left\{ \pi\left(u/K\right)\right\}
  ^{p+1}\phi\left(u,x\right)}{\{Q_{2p}\left(u/K\right)\}^{1/2}\left\{
    1+\lambda\mu\left(u\right)\right\}
}.    
\eeqn          
Properties of the Fourier transform ensure that
\beqn
\int_{-\infty}^{\infty}(2\pi\i t)^{m}\mathcal{W}(x,t)\exp(2\pi\i
tu)dt=\frac{\partial^{m}}{\partial u^{m}}\left\{ \overline{a(u,x)}\right\}.
\eeqn
Evaluating derivative of $a\left(u,x\right)$ at $u=0$ and grouping the
terms we represent
\beq   
\int_{-\infty}^{\infty}(2\pi\i t)^{m}\mathcal{W}(x,t)dt  =I_{1}+I_{2}+I_{3}, \label{eq: lemma4}
\eeq 
where 
\begin{align*}        
I_{1}=& \frac{\partial^{m}}{\partial u^{m}}\left\{\exp\left(2\pi\i x u\right)\frac{ \sinc\left(\pi u/K\right)^{2p+2}}{Q_{2p}\left(u/K\right)}\right\}_{u=0}\\
I_{2}= & \frac{\partial^{m}}{\partial u^{m}}\left[\frac{\left\{\sin\left(\pi u/K\right)\sinc\left( \pi u/K\right) \right\}^{p+1}}{Q_{2p}\left(u/K\right)}\sum_{l\neq0}\frac{\exp\left\{2\pi \i x\left(u+lK\right)\right\} }{\left\{\left(-1\right)^{l} \pi\left(u/K+l\right)\right\} ^{p+1}}\right]_{u=0}          \\ 
I_{3}=& \frac{\partial^{m}}{\partial u^{m}}\left[-\frac{\lambda\mu\left(u\right)\sinc
 \left(\pi u/K\right) ^{p+1}\phi\left(u,x\right)}{\{Q_{2p}\left(u/K\right)\}^{1/2}\left\{ 1+\lambda\mu\left(u\right)\right\} }\right]_{u=0}.  
\end{align*}
The idea is to represent each of these components as a product
of the $\sin\left(\pi u/K\right)^{n}$, $n\in\mathbb{Z}$ and some function
that is differentiable
at $0$. Then, we use that $Q_{2p}(0)=\phi(0,x)=1$, $\mu(0)=0$,
\[
\left.\frac{\partial^{m}}{\partial u^{m}} \sin\left(\pi u/K\right)^{n}\right| _{u=0}=\begin{cases}
0, & m=0,\ldots n-1,\\
n!\left(\pi/K\right)^{n}, & m=n,
\end{cases}
\]
and the Fourier series of the periodic Bernoulli polynomials
$\mathcal{B}_{p+1}(\{x\})=(-1)^{p}(p+1)!\sum_{s\neq0}\exp(-2\pi\i
sx)/\left(2\pi\i s\right)^{p+1}$.

Putting it all together and noting that $Q_{2p}(z)=\sinc(\pi
z)^{2p+2}+\sum_{l\neq 0}\sinc\{\pi(z+l)\}^{2p+2}$, we get
$
I_{1}=\left(2\pi \i x\right)^{m}$ ($m=0,\ldots, p+1$).                   
The expression for $I_{2}$ follows immediately from its representation
\[
I_{2}=\begin{cases}
0, & m=0,\ldots p,\\
-(2\pi\i/K)^{p+1}\mathcal{B}_{p+1}\left(\left\{ Kx+\frac{p+1}{2}\right\} \right), & m=p+1.
\end{cases}
\]
To find $I_{3}$, we use   $\mu\left(u\right)=(2K)^{2q}\sin(\pi u/K)^{2q}Q_{2p-2q}(u/K)/Q_{p,M}(u/K)$
\[
I_{3}=\begin{cases}
0, & m=0,\ldots2q-1,\\
-\lambda\left(2\pi\right)^{2q}\left(2q\right)!, & m=2q.       
\end{cases}
\]
To get the result for $\int_{-\infty}^{\infty}\left(t-x\right)^{m}\mathcal{W}\left(x,t\right)dt$  one needs to expand $(t-x)^m$ and use  (\ref{eq: lemma4}). 
\hfill$\square$
\subsection{Proof of Theorem \ref{theorem2} }
The asymptotic equivalent kernel on $\mathbb{R}$ can be written as
\beqn
{\cal{W}}(x,t)=\int_0^{K}\frac{\phi(u,x)\overline{\phi(u,t)}}{1+\lambda\mu(u)}du={\mathfrak{R}}\int_0^{1/2}\frac{2K\phi(Ku,x)\overline{\phi(Ku,t)}}{1+\lambda\mu(Ku)}du.
\eeqn
First consider $0\leq k_q<1$. Scaling ${\cal{W}}(x,t)$ with $c_1^{-1}K^{-1}$ leads to
\beqn
c_1{\cal{K}}(c_1x,c_1t)&=&{\mathfrak{R}}\left[\int_0^{1/2}2\phi(Ku,x/K)\overline{\phi(Ku,t/K)}du\right.\\
&-&\left.\int_0^{1/2}\frac{2\lambda\mu(Ku)
  \phi(Ku,x/K)\overline{\phi(Ku,t/K)}}{1+\lambda\mu(Ku)}du\right]={\cal{K}}_{rs}(x,t)-k_q^{2q}{\cal{K}}_1(x,t),
\eeqn
where $\mathfrak{R}$ denotes the real part of a complex number, ${\cal{K}}_{rs}(x,t)$ is the equivalent regression spline kernel on
$\mathbb{R}$ and 
\beqn
{\cal{K}}_1(x,t)&=&{\mathfrak{R}}\int_0^{1/2}\frac{2\sin(\pi
  u)^{2q}Q_{2q-2}(u)\phi(Ku,x/K)\overline{\phi(Ku,t/K)}}{\pi^{2q}Q_{2p}(u)\{1+\lambda\mu(Ku)\}}du\\
  &\leq&\frac{2^{2q}Q_{2q-2}(1/2)}{\pi^{2q}Q_{4q-2}(1/2)}{\cal{K}}_{rs}(x,t).
\eeqn
Using $Q_{lq-2}(1/2)=2\pi^{lq}(2^{lq}-1)\zeta(lq)$ for the Riemann
zeta function $\zeta(lq)=\sum_{i=1}^\infty i^{-lq}$, one can get
explicit bounds for each $q$. 

For $k_q\geq 1$ we first introduce the notation:
 $1+\lambda\mu(Ku)=\{1+\lambda(2\pi Ku)^{2q}\}\{1+r_1(u)\}$;
 $\phi(Ku,x)\overline{\phi(Ku,t)}=\exp\{2\pi\i Ku(x-t)\}\{1+r_2(x,t,u)\}$;
 $r_q(x,t,u)=\{r_2(x,t,u)-r_1(u)\}\{1+r_1(u)\}^{-1}$. 
Scaling ${\cal{W}}(x,t)$ with $c_2^{-1}\lambda^{1/(2q)}$ results in
\beqn
c_2{\cal{K}}(c_2x,c_2t)&=&\int_{-\infty}^{\infty}\frac{\exp\{2\pi\i
  u(t-x)\}}{1+(2\pi u)^{2q}}du+{\mathfrak{R}}\int_0^{k_q/2}\frac{2\exp\{2\i
  u(t-x)\}}{\pi\{1+(2u)^{2q}\}}r_q(u/k_q)du\\
&-&{\mathfrak{R}}\int_{k_q/2}^{\infty}\frac{2\exp\{2\i
  u(t-x)\}}{\pi\{1+(2u)^{2q}\}}du\\
  &=&{\cal{K}}_{ss}(x,t)+k_q^{-2q+1}{\cal{K}}_2(x,t),
\eeqn
where ${\cal{K}}_{ss}(x,t)$ is the smoothing spline kernel on
$\mathbb{R}$ and 
\beqn
\pi{\cal{K}}_2(x,t)&=&k_q^{2q-1}{\mathfrak{R}}\int_0^{k_q/2}\frac{2\exp\{2\i
  u(t-x)\}}{\pi\{1+(2u)^{2q}\}}r_q(u/k_q)du-\int_1^\infty\frac{\cos\{k_qu(t-x)\}}{\pi\{k_q^{-2q}+u^{2q}\}}du.
\eeqn 
The second component of $\pi{\cal{K}}_2(x,t)$ is obviously bounded by
$1$. Now, let us consider $r_q(u/k_q)$. First, 
$$
r_1(u/k_q)= \frac{(2u)^{2q}}{1+(2u)^{2q}}\left\{\frac{\sinc(\pi
u/k_q)^{2q}Q_{2q-2}(u/k_q)}{Q_{4q-2}(u/k_q)}-1\right\},
$$
where
$$
{Q_{2q-2}(u/k_q)}{\sinc(\pi u/k_q)^{2q}}/Q_{4q-2}(u/k_q)-1= 2\zeta(2q)(
u/k_q)^{2q}+4\zeta(4q)(u/k_q)^{4q}+\ldots
$$ 
is a positive number for any $u$. Further, 
$$r_2(x,t,u/k_q)\leq Q_{2q}^2(u/k_q)Q_{4q-2}(u/k_q)^{-1}-1=4\zeta(2q)(
u/k_q)^{2q}+8\zeta(4q)(u/k_q)^{4q}+\ldots.$$
With this,
$r_q(x,t,u/k_q)\leq 4\zeta(2q)(u/k_q)^{-2q}+\ldots$
and hence, the first term in $\pi{\cal{K}}_2(x,t)$ is also bounded for
any $k_q\geq 1$. 

Finally, ${\cal{K}}_{ss}(x,t)$ is given in \citet{Agnan96} and
${\cal{K}}_{rs}(x,t)$ is obtained from (\ref{eq:kernelW}), scaling ${\cal{W}}(x,t)$ with $K$ and setting $P_{2p}(u)=\Pi_{2p}(u)$.
\hfill$\square$
\subsection{Proof of Theorem \ref{theorem3}}
The following lemma will be used in the proof of Theorem
\ref{theorem3}.
\setcounter{lemma}{3}
\begin{lemma}
\label{lemma5} Kernel ${\cal{K}}(x,t)$, $x,t\in\mathbb{R}$ decays exponentially, i.e., there are constants $0<C<\infty$ and $0<\gamma<1$ such that
$|{\cal{K}}(x,t)|<C\gamma^{|x-t|}.$
\end{lemma}
{\sl Proof of Lemma \ref{lemma5}}\\
Since ${\cal{K}}(x,t)$ is defined as a scaled with $h(k_q)$ function
${\cal{W}}(x,t)$, from (\ref{eq:kernelW}) and (\ref{eq:bandwidth})
one finds for $k_q<1$ 
\beqn
c_{1}{\cal K}(c_{1}x,c_{1}t)=
\sum_{j=1}^{p}\sum_{l=0}^{2p}\frac{\alpha_{l}\left(\left\{ x\right\}
    ,\left\{ t\right\}
  \right)}{P'_{2p}(r_{j})}r_{j}^{\left|\left\lfloor x\right\rfloor
    -\left\lfloor t\right\rfloor +l-1 \right|+(2p-2)\mathbb{I}(\lfloor
    x\rfloor-\lfloor t\rfloor+l\leq 0)},
\eeqn
while for $k_q\geq 1$, 
\beqn
\pi c_{2}{\cal K}(\pi c_{2}x,\pi c_{2}t)=k_{q}\sum_{j=1}^{p}\sum_{l=0}^{2p}\frac{\alpha_{l}\left(\{xk_{q}\},\{tk_{q}\}\right)}{P'_{2p}(r_{j})}r_{j}^{\left|\left\lfloor xk_{q}\right\rfloor -\left\lfloor tk_{q}\right\rfloor +l-1\right|+(2p-2)\mathbb{I}(\lfloor
    xk_q\rfloor-\lfloor tk_q\rfloor+l\leq 0 )},
\eeqn
where $P_{2p}$ is given in (\ref{eq:P2p}) and $r_j=r_j(k_q)$ is a
root of  $P_{2p}$ with $\left|r_{j}\right|<1$.
If $k_q$ is a bounded constant then
$r_j=r_j(k_q)\nrightarrow\exp(-2\pi \i u),\ u\in(0,1)$ since 
$$P_{2p}\left\{ \exp(-2\pi\i u)\right\} =\exp(-2\pi\i pu)\left\{
  Q_{p,M}(u)+\left(2k_{q}/\pi\right)^{2q}\sin(\pi
  u)^{2q}Q_{2p-2q}(u)\right\} \neq0,$$ 
where the relationship between Euler--Frobenius and $Q$-polynomials has
been used. Similarly,
$r_j=r_j(k_q)\nrightarrow 0$ and $0<\gamma<1$ can be defined as 
$$\gamma=
\begin{cases}
\sup_{j,k_{q}}\left|r_{j}(k_{q})\right|,&\ k_{q}<1\\
\sup_{j,k_{q}}\left|r_{j}(k_{q})^{k_{q}}\right|,&1\leq k_{q}<\infty,
\end{cases}
$$
while 
$$
C=\sup_{k_{q},j}\frac{p(2p+1) \sup_{l,x,t}\alpha_{l}(\{x\},\{t\})}{\left|P_{2p}^{'}\{r_{j}(k_{q})\}\right|\left|r_{j}(k_q)\right|^{l+1}}
<\infty.
$$
For  $k_{q}\rightarrow\infty$ it is known from Theorem \ref{theorem2}
that 
$
\lim_{k_{q}\rightarrow\infty}{\cal K}(x,t)={\cal K}_{ss}\left\{(x-t)/\widetilde{c}_2\right\}/\widetilde{c}_2
$.
To obtain the bound on the smoothing spline kernel
${\cal{K}}_{ss}(x)$, the expression given in Theorem \ref{theorem2}
can be rewritten as
\beqn
\left|{\cal K}_{ss}(x-t)\right|&= & \left|-I_{\{{q\; is\;
      odd}\}}\frac{\exp\left(-\left|x-t\right|\right)}{2q}\right.+\sum_{j=0}^{\left\lfloor
      (q-1)/2\right\rfloor }\frac{\exp\left[-\left|x-t\right|\sin\left\{
        \pi(2j+1)/(2q)\right\} \right]}{q}\\
&&\times\left.\sin\left[\frac{\pi(2q-1)(2j+1)}{2q}-\left|x-t\right|\cos\left\{
    \frac{\pi(2j+1)}{2q}\right\} \right]
\right|\\
&\leq&\frac{q+1}{2q}\exp\left\{ -\left|x-t\right|\sin\left({\pi}{2q}\right)\right\},
\eeqn
so one can set $\gamma=\exp\left[-\sin\{\pi/(2q)\}/\widetilde{c}_2\right]\in(0,1)$, $C=(q+1)/(2q \widetilde{c}_2)<\infty$ for $k_{q}\rightarrow\infty$. 
\hfill$\square$\\\\
{\sl Proof of Theorem \ref{theorem3}}\\
Let $
\widehat{f}(x)={N}^{-1}\sum_{i=1}^NW^{[0,1]}(x,i/N)Y_i$, $\widehat{f}_{\tt{per}}(x)={N}^{-1}\sum_{i=1}^NW^{[0,1]}_{\tt{per}}(x,i/N)Y_i$.
Then, extending $f$ to the whole real line, such that it still
satisfies assumptions of the theorem, we get
\beqn
{{E}}\left\{\widehat{f}(x)\right\}&=&\int_{-\infty}^\infty{\cal{W}}(x,t)f(t)dt+R_1(x)+R_2(x)+O(N^{-1})\\
{{E}}\left\{\widehat{f}_{\tt{per}}(x)\right\}&=&\int_{-\infty}^\infty{\cal{W}}(x,t)f(t)dt+R_3(x)+O(N^{-1}),
\eeqn
where $
R_1(x)=\int_0^1\left\{ W^{[0,1]}(x,t)-{\cal{W}}(x,t)
\right\}f(t)dt$, $R_2(x)=\int_{\mathbb{R}\setminus
  [0,1]}{\cal{W}}(x,t)f(t)dt$ and $R_3(x)=\int_0^1\left\{
  W_{\tt{per}}^{[0,1]}(x,t)-{\cal{W}}_{\tt{per}}^{[0,1]}(x,t)
\right\}f(t)dt$. \\
Expanding $f(t)$ in a Taylor series around $x$ and using Lemma \ref{lemma4} results in
\beqn
\int_{-\infty}^\infty{\cal{W}}(x,t)f(t)dt&=&f(x)+ \int_{-\infty}^{\infty}\mathcal{W}\left(x,t\right)(x-t)^{2q}\frac{f^{(2q)}(\xi_{x,t})}{(2q)!}dt+O(N^{-1})\\
&= &
\frac{f^{(2q)}(x)}{(2q)!}\int_{-\infty}^{\infty}\mathcal{W}\left(x,t\right)(x-t)^{2q}dt+R_\xi(x)+O(N^{-1})\\
&=&h(k_q)^{2q}\frac{f^{(2q)}(x)}{(2q)!}\int_{-\infty}^{\infty}{\cal{K}}\left({x_h},t_h\right)\left({x_h}-t_h\right)^{2q}dt_h+R_\xi(x)+O(N^{-1}),
\eeqn
where $\xi_{x,t}$ is a point between $x$ and $t$,
$\int_{-\infty}^{\infty}{\cal{K}}\left({x_h},t\right)\left({x_h}-t\right)^{2q}dt=-C(k_q,x)$
given in the Theorem \ref{theorem3}, $x_h=x/h(k_q)$, $t_h=t/h(k_q)$ and
$$
R_\xi(x)=h(k_q)^{2q}\int_{-\infty}^{\infty}\mathcal{K}\left(x_h,t_h\right) (x_h-t_h)^{2q}\frac{f^{(2q)}(\xi_{x,t})-f^{(2q)}(x)}{h(k_q)(2q)!}dt.
$$ 
It remains to show that error terms $R_1(x)$, $R_2(x)$ are
negligible for $x\in{\cal{I}}_q$ and $R_3(x)$,
  $R_\xi(x)$ are uniformly negligible. \\Using techniques similar to \citet{HuangStudden93},
\beqn
R_\xi(x)&= & h(k_q)^{2q}\sum_{l=-\infty}^{\infty}\int_{x+(l-1)h(k_q)}^{x+lh(k_q)}\mathcal{K}\left(x_h,{t}_h\right) (x_h-t_h)^{2q}\frac{f^{(2q)}(\xi_{x ,t})-f^{(2q)}(x)}{h(k_q)(2q)!}dt\\
&\leq&
h(k_q)^{2q+\alpha}CL\sum_{l=-\infty}^{\infty}\int_{x+(l-1)h(k_q)}^{x+lh(k_q)}\gamma^{|x_h-t_h|}\frac{\left|x_h-t_h\right|^{2q+\alpha}}{h(k_q)(2q)!}dt\\
&\leq&h(k_q)^{2q+\alpha}\frac{2CL}{(2q)!}\sum_{l=1}^\infty\gamma^{l-1}l^{2q+\alpha}=o\left\{h(k_q)^{2q}\right\},
\eeqn
where the exponential bound on the kernel from Lemma \ref{lemma5}
together with 
the H\"older continuity of $f^{(2q)}$ have been used.

To see that
$R_3(x)=o\left[\{h(k_q)N\}^{-1}\right]$ for any $x$,
use the definitions of both kernels to get
\beqn
|R_3(x)|&\leq&\|f\|_\infty\int_0^1\left|\sum_{i=1}^K\frac{\phi_i(x)\overline{\phi_i(t)}}{1+\lambda
\mu(i)}\frac{Q_{p,M}(i/K)-Q_{2p}(i/K)}{Q_{p,M}(i/K)\{1+\lambda\nu(i)\}}\right|dt=O\left[\{h(k_q)N\}^{-2q}\right],
\eeqn 
since $Q_{p,M}(i/K)=Q_{2p}(i/K)+\sin(\pi i/K)^{2q}M^{-2q}$ and both $Q_{2p}(i/K),Q_{p,M}(i/K)\in(0,1]$ for any $i=1,\ldots,K$.
Next,
\beqn
|R_2(x)|&=&\left|\int_{\mathbb{R}\setminus[0,1]}{\cal{W}}(x,t)f(t)dt\right|\leq\|f\|_\infty\int_{\mathbb{R}\setminus[0,1]}h(k_q)^{-1}\gamma^{|x-t|h(k_q)^{-1}}dt\\
&=&\|f\|_\infty\left( \gamma^{\frac{x}{h(k_q)}}+\gamma^{\frac{1-x}{h(k_q)}} \right)/\log(1/\gamma)=o\left\{h(k_q)^{2q}\right\}
\eeqn
as long as $x\in{\cal{I}}_q$ with
$\delta_q>2q\log_\gamma(1/e)$. Finally, for $x\in{\cal{I}}_q$,
\beqn
R_1(x)&=&\int_0^1\left\{ W^{[0,1]}(x,t)-{\cal{W}}(x,t)
\right\}f(t)dt\\
&=&\int_0^1\left\{W^{[0,1]}(x,t)-W^{[0,1]}_{\tt{per}}(x,t)\right\}f(t)dt+R_2(x)+R_3(x)\\
&=&O(N^{-1})+o\left\{h(k_q)^{2q}\right\}+O\left[\{h(k_q)N\}^{-2q}\right],
\eeqn
since the difference between
 projection of function $f$ onto general and periodic spline
spaces is zero by definition for $x\in[2q/K,1-2q/K]\supset{\cal{I}}_q$, see Chapter 8.1
in \citet{Schumaker07}.

Now, the variance of $\widehat{f}_{\tt{per}}(x)$ is given by
\[
\mbox{var}\left\{\widehat{f}_{\tt{per}}(x)\right\}=\frac{\sigma^{2}}{N^{2}}\sum_{i=1}^{N}W^{[0,1]}_{\tt{per}}(x,i/N)^2=\frac{\sigma^{2}}{N}\int_{0}^{1}{\cal{W}}^{[0,1]}_{\tt{per}}(x,t)^2dt+\frac{\sigma^2}{N}R_4(x)+O(N^{-2}).
\]
for $R_4(x)=\int_0^1\left\{{{W}}^{[0,1]}_{\tt{per}}(x,t)^2-{\cal{W}}^{[0,1]}_{\tt{per}}(x,t)^2\right\}dt$.
Let us define ${\tt{K}}_{\tt{per}}(x,t)$ via 
$$h(k_q)^{-1}{\tt{K}}_{\tt{per}} (x_h,t_h)={\cal{W}}^{[0,1]}_{\tt{per}}(x,t)=\sum_{l=-\infty}^\infty{\cal{W}}(x,t+l)=h(k_q)^{-1}\sum_{l=-\infty}^\infty{\cal{K}}(x_h,t_h+l_h),$$
for $l_h=l/h(k_q)$. Then, using periodicity of ${\cal{W}}_{\tt{per}}^{[0,1]}(x,t)$
\beqn
\int_0^1{\cal{W}}^{[0,1]}_{\tt{per}}(x,t)^2dt&=&\int_{x-1/2}^{x+1/2}{\cal{W}}^{[0,1]}_{\tt{per}}(x,t)^2dt=\frac{1}{h(k_q)^2}\int_{x-1/2}^{x+1/2}{\tt{K}}_{\tt{per}}(x_h,t_h)^2dt\\
&=&\frac{1}{h(k_q)}\left\{\int_{-\infty}^\infty{\cal{K}}(x_h,t)^2dt+ R_k(x)\right\},
\eeqn
for 
\beqn
h(k_q)
R_k(x)&=&\int_{x-1/2}^{x+1/2}{\tt{K}}_{\tt{per}}(x_h,t_h)^2dt-\int_{-\infty}^\infty{\cal{K}}(x_h,t_h)^2dt\\
&=&\int_{x-1/2}^{x+1/2}\left\{{\tt{K}}_{\tt{per}}(x_h,t_h)^2-{\cal{K}}(x_h,t_h)^2\right\}dt\\
&-&\int_{-\infty}^{x-1/2}{\cal{K}}_{\tt{per}}(x_h,t_h)^2dt-\int_{x+1/2}^\infty{\cal{K}}(x_h,t_h)^2dt.
\eeqn
Now, we can make use of
${\tt{K}}_{\tt{per}} (x,t)=\sum_{l=-\infty}^\infty{\cal{K}}(x,t+l)$ and of the
exponential decay of ${\cal{K}}(x,t)$ found in Lemma
\ref{lemma5} to bound terms in $h(k_q) R_k(x)$. That is,
\beqn
\int_{x-1/2}^{x+1/2}\left\{{\tt{K}}_{\tt{per}}(x_h,t_h)^2-{\cal{K}}(x_h,t_h)^2\right\}dt&=&\int_{x-1/2}^{x+1/2}\sum_{l\neq 0}{\cal{K}}(x_h,t_h+l_h) \left\{ \sum_{l\neq
    0}{\cal{K}}(x_h,t_h+l_h)+2{\cal{K}}(x_h,t_h) \right\}dt\\
&\leq&C^2
\int_{x-1/2}^{x+1/2}\sum_{l\neq 0}\gamma^{|x_h-t_h-l_h|}\left(\sum_{l\neq
    0}\gamma^{|x_h-t_h-l_h|}  +2\gamma^{|x_h-t_h|}\right)dt\\
&\leq& h(k_q)\frac{C^2 \gamma^{1/h(k_q)}\left\{
    4+2h(k_q)^{-1}\log(\gamma^{-1})\right\}}{\left\{\gamma^{1/h(k_q)}-1\right\}^2\log(\gamma^{-1})},
\eeqn
where 
$
\sum_{l\neq 0}\gamma^{|x_h-t_h-l_h|}=\left(\gamma^{t_h-x_h}+\gamma^{x_h-t_h}\right)\gamma^{1/h(k_q)}/\left\{1-\gamma^{1/h(k_q)}\right\},
$
for $t\in[x-1/2,x+1/2]$ has been used. Also, 
\beqn
\int_{-\infty}^{x-1/2}{\cal{K}}(x_h,t_h)^2dt+\int_{x+1/2}^\infty{\cal{K}}(x_h,t_h)^2dt&\leq&
C^2\left\{\int_{-\infty}^{x-1/2}\gamma^{2(x_h-t_h)}dt+\int_{x+1/2}^\infty\gamma^{2(t_h-x_h)}dt\right\}\\
&=&h(k_q)\frac{C^2\gamma^{{1}/{h(k_q)}}}{\log(\gamma^{-1})}.
\eeqn
In a similar fashion one finds 
$\int_{-\infty}^{\infty}{\cal{K}}(x_h,t)^2dt\leq
{C^2}/{\log(\gamma^{-1})}$. Putting it all together gives 
\beqn
\left|R_k(x)\right|&\leq&
\frac{C^2\gamma^{1/h(k_q)}}{\log(\gamma^{-1})}\left[1+\frac{4+2h(k_q)^{-1}\log(\gamma^{-1})}{\left\{\gamma^{1/h(k_q)}-1
    \right\}^2} \right]=O\left\{h(k_q)^{-1}\gamma^{1/h(k_q)}\right\}=o(1).
\eeqn
The proof for $\mbox{var}\left\{\widehat{f}(x)\right\}$ follows from
\beqn
\mbox{var}\left\{\widehat{f}(x)\right\}&=&\frac{\sigma^2}{N}\sum_{i=1}^NW^{[0,1]}(x,i/N)^2=\frac{\sigma^2}{N}\int_{-\infty}^{\infty}{\cal{W}}(x,t)^2dt+\frac{\sigma^2}{N}\left\{R_5(x)+R_6(x)\right\}+O\left(N^{-2}\right)\\
&=&\frac{\sigma^2}{Nh(k_q)}\int_{-\infty}^\infty{\cal{K}}\{x/h(k_q),t\}^2dt+\frac{\sigma^2}{N}\left\{R_5(x)+R_6(x)\right\}+O\left(N^{-2}\right),
\eeqn
where
$R_5(x)=\int_0^1\left\{W^{[0,1]}(x,t)^2-{\cal{W}}(x,t)^2\right\}dt$
and $R_6(x)=\int_{\mathbb{R}\setminus[0,1]}{\cal{W}}(x,t)^2dt$. The
proof that $R_4(x) $ is uniformly negligible and $R_5(x)$, $R_6(x)$
are negligible for $x\in {\cal{I}}_q$ follows exactly the same lines
as that for $R_3(x)$, $R_1(x)$ and $R_2(x)$, respectively.
\hfill$\square$
\end{appendix}
\biblist

\end{document}